\documentclass{amsart}
\usepackage[latin9]{inputenc}
\usepackage{geometry}
\usepackage{amstext}
\usepackage{amsthm}
\usepackage{amssymb}
\usepackage{graphicx}
\usepackage{esint}

\makeatletter

\DeclareTextSymbolDefault{\textquotedbl}{T1}

\numberwithin{equation}{section}
\numberwithin{figure}{section}
\theoremstyle{plain}
\newtheorem{thm}{\protect\theoremname}
\theoremstyle{plain}
\newtheorem{lem}[thm]{\protect\lemmaname}

\makeatother

\providecommand{\lemmaname}{Lemma}
\providecommand{\theoremname}{Theorem}

\begin{document}
\title{Zeros of a binomial combination of Chebyshev polynomials}
\author{Summer Al Hamdani \and Khang Tran}
\address{Department of Mathematics \\
California State University\\
5245 North Backer Avenue M/S PB108 Fresno, CA 93740}
\begin{abstract}
For $0<\alpha<1$, we study the zeros of the sequence of polynomials
$\left\{ P_{m}(z)\right\} _{m=0}^{\infty}$ generated by the reciprocal
of $(1-t)^{\alpha}(1-2zt+t^{2})$, expanded as a power series in $t$.
Equivalently, this sequence is obtained from a linear combination
of Chebyshev polynomials whose coefficients have a binomial form.
We show that the number of zeros of $P_{m}(z)$ outside the interval
$(-1,1)$ is bounded by a constant independent of $m$. 
\end{abstract}

\maketitle

\section{Introduction}

\begin{figure}
\begin{centering}
\includegraphics[scale=0.5]{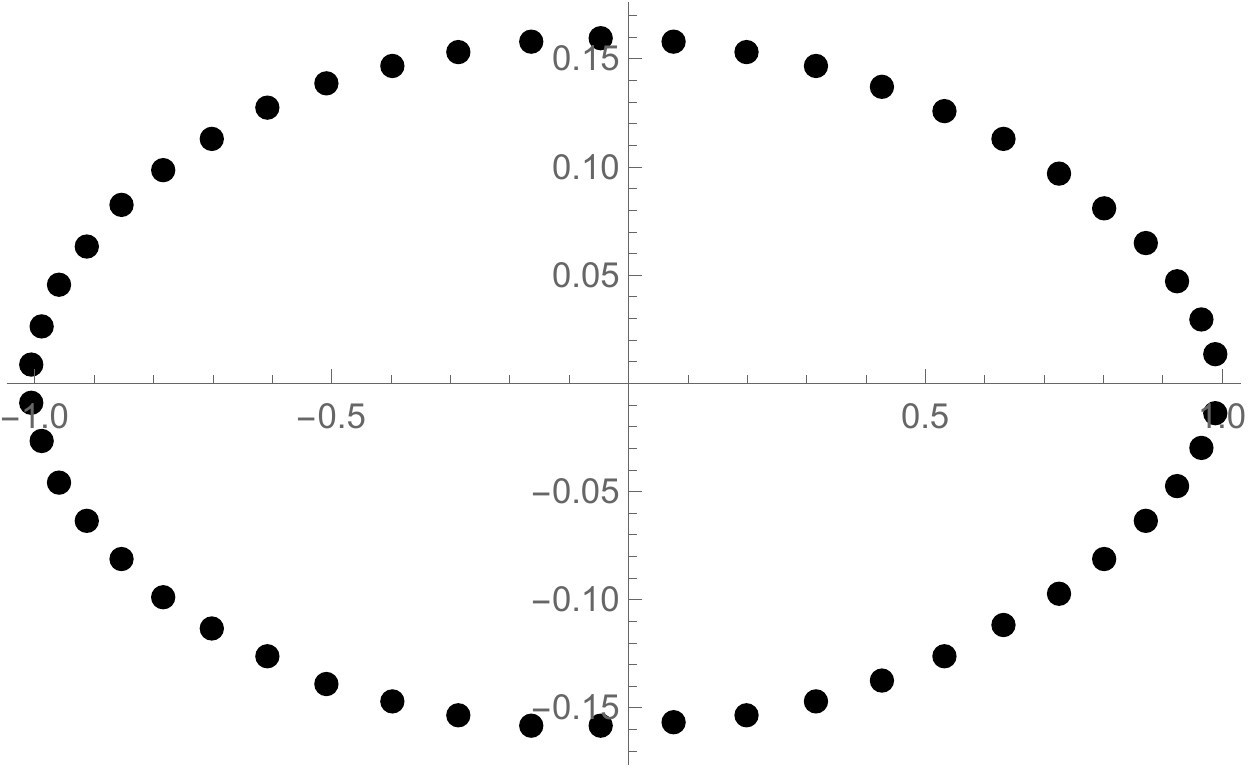}
\par\end{centering}
\caption{\label{fig:example}Zeros of $P_{50}(z)$ for $\alpha=3$}
\end{figure}

The sequences of Chebyshev polynomials of the first and second kinds,
generated respectively by 
\[
\sum_{m=0}^{\infty}T_{m}(z)t^{m}=\frac{1-tz}{1-2tz+t^{2}}
\]
and 
\[
\sum_{m=0}^{\infty}U_{m}(z)t^{m}=\frac{1}{1-2tz+t^{2}},
\]
are two important sequences of orthogonal polynomials in mathematics.
This orthogonality implies that the zeros of all polynomials in these
sequences lie on the support of the weight function which is the real
interval $(-1,1)$. For each $n\in\mathbb{N}$, the sum of the first
$n$ Chebyshev polynomials of the first kind connects with the famous
Dirichlet kernel, $D_{n}(\theta)$, by 
\[
D_{n}(\theta)=\frac{1}{2\pi}\left(2\sum_{m=0}^{n}T_{m}(z)-1\right)=\frac{\sin((n+1/2)\theta)}{2\pi\sin(\theta/2)}
\]
where $z=\cos\theta$. One can use this relation to show that the
zeros of the partial sum $\sum_{m=0}^{n}T_{m}(z)$ still lie on the
interval $(-1,1)$. With a similar identity, we can show that the
same conclusion holds for the partial sum of Chebyshev polynomials
of the second kind, $V_{n}(z)=\sum_{m=0}^{n}U_{m}(z)$. For a study
of the distribution of zeros of the sum of a more general sequence
of polynomials satisfying a three-term recurrence, see \cite{tz}.
The generating for this sequence of partial sums $\left\{ V_{m}(z)\right\} _{m=0}^{\infty}$
is given by 
\begin{align}
\sum_{n=0}^{\infty}V_{n}(z)t^{n}= & \sum_{n=0}^{\infty}\sum_{m=0}^{n}U_{m}(z)t^{n}\nonumber \\
= & \sum_{m=0}^{\infty}U_{m}(z)t^{m}\sum_{n=m}^{\infty}t^{n-m}\nonumber \\
= & \frac{1}{(1-t)(1-2tz+t^{2})}.\label{eq:sumgenfunc}
\end{align}
We can look at the sequence of partial sums of $V_{m}(z)$ (i.e. the
sequence $\left\{ \sum_{m=0}^{n}V_{m}(z)\right\} _{n=0}^{\infty}$)
or more generally if we iterate the sequence of partial sums $\alpha$
times then the generating function of the resulting sequence, denoted
by $\left\{ P_{m}(z)\right\} _{m=0}^{\infty}$, is 
\[
\sum_{m=0}^{\infty}P_{m}(z)t^{m}=\frac{1}{(1-t)^{\alpha}(1-2tz+t^{2})}.
\]
However the zeros of polynomials $P_{m}(z)$ in this sequence for
$\alpha\in\mathbb{N}$ may not lie on the interval $(-1,1)$ (see
Figure \ref{fig:example}). Surprisingly, the zeros of $P_{m}(z)$
are more likely to lie on this interval when $0<\alpha<1$. The main
goal of this paper is to prove the theorem below.
\begin{thm}
\label{thm:mainthm}For $0<\alpha<1$, $\alpha\in\mathbb{R}$, let
$\left\{ P_{m}(z)\right\} _{m=0}^{\infty}$ be the sequence of polynomials
generated by 
\[
\sum_{m=0}^{\infty}P_{m}(z)t^{m}=\frac{1}{(1-t)^{\alpha}(1-2zt+t^{2})}.
\]
There is a constant $C$ (independent of $m$) such that the number
of zeros of $P_{m}(z)$ outside $(-1,1)$ is at most $C$ for all
$m\in\mathbb{N}$. 
\end{thm}

One can apply a procedure similar to \eqref{eq:sumgenfunc} to write
$P_{m}(z)$ as the following binomial combination of Chebyshev polynomials
of the second kind. 
\begin{equation}
P_{m}(z)=\sum_{k=0}^{m}\binom{\alpha+m-k-1}{m-k}U_{k}(z).\label{eq:binomcomb}
\end{equation}
Although we do not work with such binomial combination of Chebyshev
polynomials of the first, the method in this paper should work for
that combination. We focus on binomial combinations since binomial
coefficients play an important role in the research program on stability
preserving linear operators on circular domains initiated by P\'{o}lya
and Schur and completed by J. Borcea and P. Br\"{a}nd\'{e}n \cite{bb}.
 For a study of zeros of other linear combinations of Chebyshev polynomials,
see \cite{stankov}. The locations of zeros of polynomials strongly
relate to the discriminants and resultants of those polynomials. Given
that the discriminants and resultants of Chebyshev polynomials are
known (see \cite{gi,rivlin,szego}), there have been various studies
extending these concepts to various linear combinations of Chebshev
polynomials (for example, see \cite{dk,tran}). 

\section{The Cauchy integral formula and a proof of the theorem}

\begin{figure}
\begin{centering}
\includegraphics[scale=0.3]{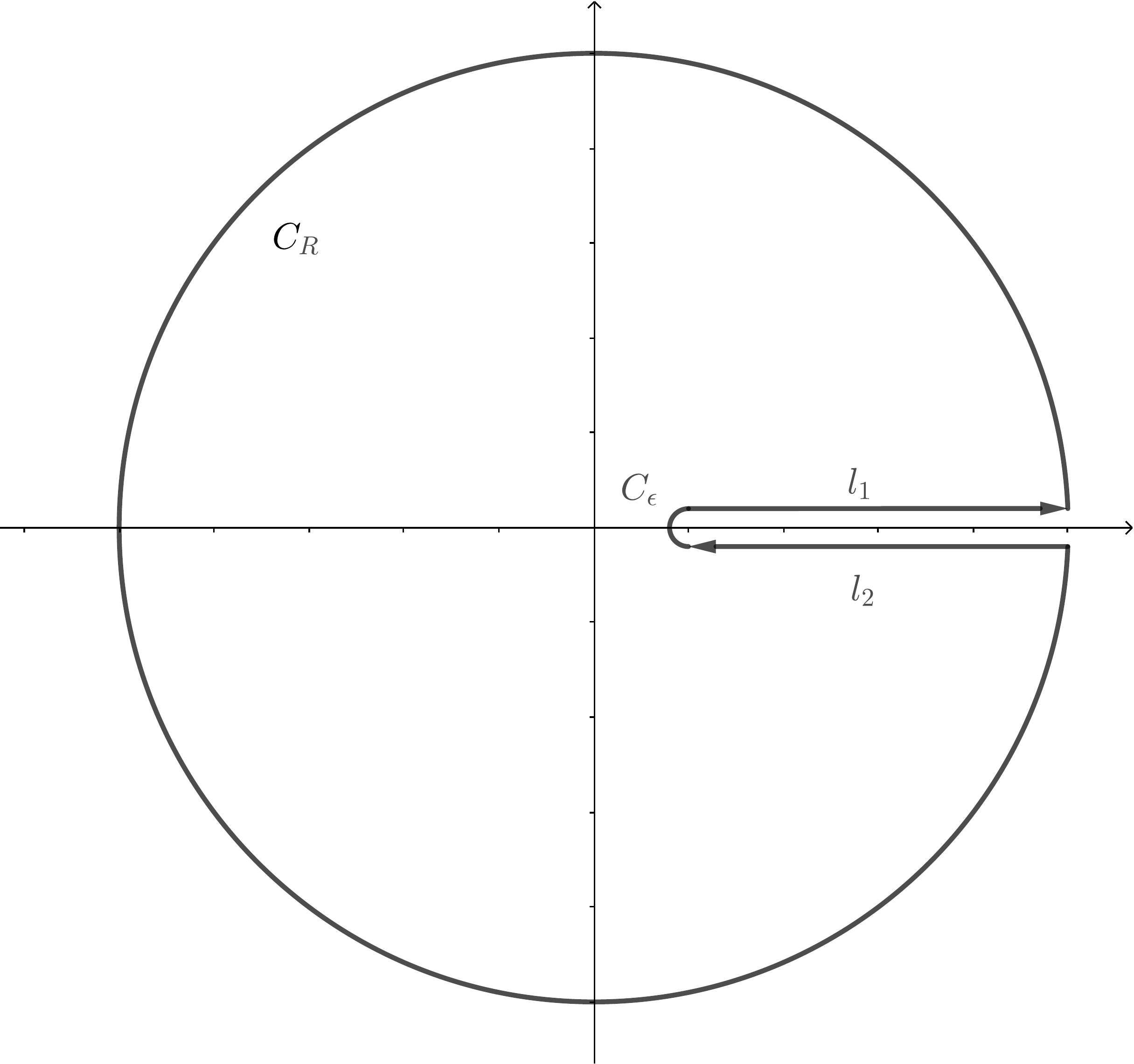}
\par\end{centering}
\caption{\label{fig:contour}The contour $\gamma$}
\end{figure}

We start our proof by defining the function $z(\theta)=\cos\theta$
on the interval $\theta\in(0,\pi)$ and note that, for each $\theta\in(0,\pi)$,
the two zeros in $t$ of $1+2z(\theta)t+t^{2}$ are $e^{\pm i\theta}$.
With the Cauchy differentiation formula and the principal cut for
$(1-t)^{\alpha}$, we conclude that for each $\theta\in(0,\pi)$
\begin{align}
P_{m}(z(\theta)) & =\frac{1}{2\pi i}\ointctrclockwise_{|t|=\delta}\frac{dt}{(1-t)^{\alpha}(1-2z(\theta)t+t^{2})t^{m+1}}\nonumber \\
 & =\frac{1}{2\pi i}\ointctrclockwise_{|t|=\delta}\frac{dt}{(1-t)^{\alpha}(t-e^{i\theta})(t-e^{-i\theta})t^{m+1}}\label{eq:cauchyint}
\end{align}
for some small $\delta>0$. For each large $R$ and small $\epsilon$,
let $\gamma$ be the counterclockwise loop formed by the union of
$C_{R}$ (a portion of the circle radius large $R$), and a small
semicircle $C_{\epsilon}$ with center $1$ and radius $\epsilon$,
and the two line segments $l_{1}$ and $l_{2}$ (see Figure \ref{fig:contour}).

We note that for each $\theta\in(0,\pi)$
\begin{align*}
 & \left|\int_{C_{R}}\frac{dt}{(1-t)^{\alpha}(t-e^{i\theta})(t-e^{-i\theta})t^{m+1}}\right|\\
\le & \int_{C_{R}}\frac{|dt|}{(R-1)^{\alpha}(R-1)(R-1)R^{m+1}}\\
\rightarrow & 0
\end{align*}
as $R\rightarrow\infty$ and 
\begin{align*}
 & \left|\int_{C_{\epsilon}}\frac{dt}{(1-t)^{\alpha}(t-e^{i\theta})(t-e^{-i\theta})t^{m+1}}\right|\\
\le & \int_{_{\pi/2}}^{3\pi/2}\frac{\epsilon d\phi}{\epsilon^{\alpha}|1+\epsilon e^{i\phi}-e^{i\theta}||1+\epsilon e^{i\phi}-e^{-i\theta}||1+\epsilon e^{i\phi}|^{m+1}}\\
\rightarrow & 0
\end{align*}
as $\epsilon\rightarrow0$ since $\alpha<1$. 

With the paramatrization $t=1+x+i\epsilon$ for $l_{1}$, we note
that for each $m$ and $\theta$ the function in $x$ 
\begin{equation}
\frac{1}{(1-t)^{\alpha}(t-e^{i\theta})(t-e^{-i\theta})t^{m+1}}\label{eq:integrand}
\end{equation}
converges pointwise to 
\[
\frac{1}{x^{\alpha}e^{-i\pi\alpha}(1+x-e^{i\theta})(1+x-e^{-i\theta})(1+x)^{m+1}}
\]
on $x\in(0,\infty)$ as $R\rightarrow\infty$ and $\epsilon\rightarrow0$.
With the note that 
\[
|1+x-e^{i\theta}||1+x-e^{-i\theta}|=(1+x)^{2}-2\cos\theta(1+x)+1
\]
is an increasing function in $x$ on $(0,\infty)$, we conclude 
\[
|1+x-e^{i\theta}||1+x-e^{-i\theta}|\ge2-2\cos\theta.
\]
Consequently, we can dominate \eqref{eq:integrand} by 
\[
\frac{1}{|1-t|^{\alpha}|t-e^{i\theta}||t-e^{-i\theta}||t|^{m+1}}\le\frac{1}{x^{\alpha}(2-2\cos\theta)(1+x)^{m+1}}
\]
where the right side is integrable on $(0,\infty)$ since $0<\alpha<1$.
Thus by the dominated convergence theorem, we conclude that as $R\rightarrow\infty$
and $\epsilon\rightarrow0$
\begin{align*}
\int_{l_{1}}\frac{dt}{(1-t)^{\alpha}(t-e^{i\theta})(t-e^{-i\theta})t^{m+1}} & \rightarrow\int_{0}^{\infty}\frac{dx}{x^{\alpha}e^{-i\pi\alpha}(1+x-e^{i\theta})(1+x-e^{-i\theta})(1+x)^{m+1}}
\end{align*}
and similarly for $l_{2}$
\[
\int_{l_{2}}\frac{dt}{(1-t)^{\alpha}(t-e^{i\theta})(1-e^{-i\theta})t^{m+1}}\rightarrow-\int_{0}^{\infty}\frac{dx}{x^{\alpha}e^{i\pi\alpha}(1+x-e^{i\theta})(1+x-e^{-i\theta})(1+x)^{m+1}}.
\]

Let $f(t,\theta)$ be the integrand of \eqref{eq:cauchyint} and note
that as a function in $t$, this function has three poles at $0$
and $e^{\pm i\theta}.$ From
\begin{align*}
 & \frac{1}{2\pi i}\int_{C_{R}}f(t,\theta)dt+\frac{1}{2\pi i}\int_{C_{\epsilon}}f(t,\theta)dt+\frac{1}{2\pi i}\int_{l_{1}}f(t,\theta)dt+\frac{1}{2\pi i}\int_{l_{2}}f(t,\theta)dt\\
= & P_{m}(z(\theta))+\ointctrclockwise_{|t-e^{i\theta}|=\delta}f(t,\theta)dt+\ointctrclockwise_{|t-e^{-i\theta}|=\delta}f(t,\theta)dt,
\end{align*}
we let we let $R\rightarrow\infty$ and $\epsilon\rightarrow0$ to
obtain 
\begin{align}
P_{m}(z(\theta)) & =\frac{1}{\pi}\Im\int_{0}^{\infty}\frac{dx}{x^{\alpha}e^{-i\pi\alpha}(1+x-e^{i\theta})(1+x-e^{-i\theta})(1+x)^{m+1}}\nonumber \\
 & -\frac{1}{2\pi i}\ointctrclockwise_{|t-e^{-i\theta}|=\delta}f(t,\theta)dt-\frac{1}{2\pi i}\ointctrclockwise_{|t-e^{i\theta}|=\delta}f(t,\theta)dt.\label{eq:PmCauchy}
\end{align}
Since $e^{\pm i\theta}$ are two simple poles of $f(t,\theta)$, the
Cauchy integral formula gives
\begin{align*}
 & -\frac{1}{2\pi i}\ointctrclockwise_{|t-e^{-i\theta}|=\delta}f(t,\theta)dt-\frac{1}{2\pi i}\ointctrclockwise_{|t-e^{i\theta}|=\delta}f(t,\theta)dt\\
= & -\frac{1}{(1-e^{-i\theta})^{\alpha}(e^{-i\theta}-e^{i\theta})e^{-i(m+1)\theta}}-\frac{1}{(1-e^{i\theta})^{\alpha}(e^{i\theta}-e^{-i\theta})e^{i(m+1)\theta}}.
\end{align*}
After forming a common denominator and applying $(1-e^{-i\theta})(1+e^{i\theta})=2-2\cos\theta$,
the last expression becomes
\[
\frac{\left(1-e^{i\theta}\right)^{\alpha}e^{i(m+1)\theta}-\left(1-e^{-i\theta}\right)^{\alpha}e^{-i(m+1)\theta}}{2i\sin\theta(2-2\cos\theta)^{\alpha}}.
\]
With the identities 
\begin{align*}
(1-e^{i\theta})^{\alpha}e^{i(m+1)\theta} & =(e^{-i\theta/2}-e^{i\theta/2})^{\alpha}e^{i(m+1)\theta+i\alpha\theta/2}\\
 & =2^{\alpha}\sin^{\alpha}(\theta/2)e^{i(m+1)\theta+i\alpha\theta/2-i\alpha\pi/2}
\end{align*}
and 
\begin{align*}
(1-e^{-i\theta})^{\alpha}e^{-i(m+1)\theta} & =(e^{i\theta/2}-e^{-i\theta/2})^{\alpha}e^{-i(m+1)\theta-i\alpha\theta/2}\\
 & =2^{\alpha}\sin^{\alpha}(\theta/2)e^{-i(m+1)\theta-i\alpha\theta/2+i\alpha\pi/2},
\end{align*}
we write this expression as
\[
\frac{\sin^{\alpha}(\theta/2)}{\sin\theta(1-\cos\theta)^{\alpha}}\sin\left((m+1)\theta+\frac{\alpha(\theta-\pi)}{2}\right).
\]
We plug this quantity to the last two terms of the right side of \eqref{eq:PmCauchy}
and conclude for $\theta\in(0,\pi)$
\begin{align}
P_{m}(z(\theta)) & =\frac{1}{\pi}\Im\int_{0}^{\infty}\frac{dx}{x^{\alpha}e^{-i\pi\alpha}((1+x)^{2}-2(1+x)\cos\theta+1)(1+x)^{m+1}}\nonumber \\
 & +\frac{\sin^{\alpha}(\theta/2)}{\sin\theta(1-\cos\theta)^{\alpha}}\sin\left((m+1)\theta+\frac{\alpha(\theta-\pi)}{2}\right).\label{eq:Pmform}
\end{align}
In the next section, we will prove the key lemma below.
\begin{lem}
\label{lem:keylemma}There are constants $K$ and $M$ such that 
\begin{equation}
\frac{1}{\pi}\int_{0}^{\infty}\frac{dx}{x^{\alpha}((1+x)^{2}-2(1+x)\cos\theta+1)(1+x)^{m+1}}<\frac{\sin^{\alpha}(\theta/2)}{\sin\theta(1-\cos\theta)^{\alpha}}\label{eq:keyineq}
\end{equation}
for all $\theta\in(K/m,\pi)$ and all $m\ge M$. 
\end{lem}

Assume the lemma above, we now provide a proof of Theorem \ref{thm:mainthm}.
Consider all the angles $\theta_{h}\in(K/m,\pi)$ such that 
\[
\sin\left((m+1)\theta_{h}+\frac{\alpha(\theta_{h}-\pi)}{2}\right)=\pm1.
\]
That is 
\[
\theta_{h}=\frac{h\pi+(\alpha+1)\pi/2}{m+1+\alpha/2}.
\]
The condition $\theta_{h}\in(K/m,\pi)$ implies
\[
\frac{K}{m}<\frac{h\pi+(\alpha+1)\pi/2}{m+1+\alpha/2}<\pi
\]
or equivalently 
\begin{equation}
\frac{m+1+\alpha/2}{m\pi}K-\frac{\alpha+1}{2}<h<m+\frac{1}{2}.\label{eq:thetahrange}
\end{equation}
Lemma \ref{lem:keylemma} and \eqref{eq:Pmform} imply that for all
$m\ge M$, the sign of $P_{m}(z(\theta_{h}))$ is the same as the
sign of 
\[
\sin\left((m+1)\theta_{h}+\frac{\alpha(\theta_{h}-\pi)}{2}\right)
\]
which is $(-1)^{h}$. By the Intermediate Value Theorem, for each
$h$ the interval $(\theta_{h},\theta_{h+1})$ contains at least a
zero of $P_{m}(z(\theta))$. From \eqref{eq:thetahrange}, we conclude
that there is constant $C_{1}$ such that $P_{m}(z(\theta))$ has
at least $m-C_{1}$ zeros in $\theta$ on $(K/m,\pi)$ for $m\ge M$
and each of these zeros gives a zero of $P_{m}(z)$ on $(-1,1)$ by
the map $z(\theta)=\cos\theta$. Let $C=\max(C_{1},M)$. Since the
degree of $P_{m}(z)$ is $m$ by \eqref{eq:binomcomb}, the number
of zeros of $P_{m}(z)$ outside $(-1,1)$ is at most $C$ for all
$m\in\mathbb{N}$. 

\section{Proof of Lemma \ref{lem:keylemma}}

In this section, we will prove Lemma \ref{lem:keylemma}. With the
substitution $1+x=e^{u}$, we write the integral on the left of \eqref{eq:keyineq}
as 
\begin{align*}
\frac{1}{\pi}\int_{0}^{\infty}\frac{dx}{x^{\alpha}(1+x)^{m+1}((1+x)^{2}-2(1+x)\cos\theta+1)} & =\frac{1}{\pi}\int_{0}^{\infty}\frac{e^{-mu}du}{(e^{u}-1)^{\alpha}(e^{2u}-2e^{u}\cos\theta+1)}\\
 & =\frac{1}{\pi}\int_{0}^{\infty}\frac{e^{-(m+\alpha)u}du}{(1-e^{-u})^{\alpha}(e^{2u}-2e^{u}\cos\theta+1)}\\
 & =\frac{1}{\pi}\int_{0}^{\infty}e^{-(m+\alpha)u}u^{-\alpha}g(u)du,
\end{align*}
where 
\begin{equation}
g(u)=\dfrac{f(u)}{e^{2u}-2e^{u}\cos\theta+1}\qquad\text{and}\qquad f(u)=\left(\dfrac{u}{1-e^{-u}}\right)^{\alpha}.\label{eq:gfdef}
\end{equation}
Note that if $\theta$ is fixed, one can apply Watson's lemma to find
an asymptotic formula (nonuniform in $\theta$) for the integral above.
Since the range of $\theta$ of interest in Lemma \ref{lem:keylemma}
depends on $m$, we cannot apply Watson's lemma directly. The book
\cite{temme} provides many uniform asymptotic formulas for various
other important integrals. Here we apply an elementary approach based
on the proof of Watson's lemma to prove Lemma \ref{lem:keylemma}.

To proceed, we split the integral above as 
\begin{equation}
\frac{1}{\pi}\int_{0}^{1/\sqrt{m}}e^{-(m+\alpha)u}u^{-\alpha}g(u)du+\frac{1}{\pi}\int_{1/\sqrt{m}}^{\infty}e^{-(m+\alpha)u}u^{-\alpha}g(u)du.\label{eq:splitrange}
\end{equation}
Note that $e^{2u}-2e^{u}\cos\theta+1$ is a decreasing function in
$u$ on $(0,\infty)$ and consequently 
\begin{align*}
e^{2u}-2e^{u}\cos\theta+1 & \ge2-2\cos\theta\\
 & \gg\theta^{2}.
\end{align*}
With this inequality and $f(u)=\mathcal{O}(u^{\alpha})$ (the Big-Oh
constant here and throughout the paper is independent of $\theta$),
we obtain the following bound for the second integral of \eqref{eq:splitrange}

\begin{align}
\left|\frac{1}{\pi}\int_{1/\sqrt{m}}^{\infty}e^{-(m+\alpha)u}u^{-\alpha}g(u)du\right| & =\mathcal{O}\left(\frac{1}{\theta^{2}}\int_{1/\sqrt{m}}^{\infty}e^{-(m+\alpha)u}du\right)\nonumber \\
 & =\mathcal{O}\left(\frac{e^{-\sqrt{m}}}{m\theta^{2}}\right).\label{eq:secondint}
\end{align}

Next we will bound the first integral of \eqref{eq:splitrange}. With
the note that 
\[
\lim_{u\rightarrow0}f(u)=\lim_{u\rightarrow0}\left(\dfrac{u}{1-e^{-u}}\right)^{\alpha}=1
\]
we can define 
\begin{equation}
g(0):=\lim_{u\rightarrow0}g(u)=\frac{1}{2-2\cos\theta}\label{eq:gat0}
\end{equation}
 and conclude $g(u)=g(0)+g'(v)u$ for some $v\in(0,u)$ where from
\eqref{eq:gfdef}
\begin{align*}
g'(v) & =\frac{(e^{2v}-2e^{v}\cos\theta+1)f'(v)-f(v)(2e^{2v}-2e^{v}\cos\theta)}{\left(e^{2v}-2e^{v}\cos\theta+1\right)^{2}}\\
 & \leq\frac{(e^{2v}-2e^{v}\cos\theta+1)f'(v)}{\left(e^{2u}-2e^{u}\cos\theta+1\right)^{2}}\\
 & =\frac{f'(v)}{e^{2v}-2e^{v}\cos\theta+1}.
\end{align*}
From the definition of $f$ in \eqref{eq:gfdef}, we have 
\[
f'(u)=\alpha\left(\dfrac{u}{1-e^{-u}}\right)^{\alpha-1}\frac{1-e^{-u}(1+u)}{\left(1-e^{-u}\right)^{2}}.
\]
The identity above implies
\[
\lim_{u\rightarrow0}f'(u)=\frac{\alpha}{2}\qquad\text{and}\qquad\lim_{u\rightarrow\infty}f'(u)=0,
\]
from which we conclude $f'(u)$ is bounded on $(0,\infty)$ and consequently
\[
g'(v)=\mathcal{O}\left(\frac{1}{e^{2v}-2e^{v}\cos\theta+1}\right)=\mathcal{O}\left(\frac{1}{\theta^{2}}\right).
\]
Together with \eqref{eq:gat0}, we conclude for $0<u<1/\sqrt{m}$
\[
g(u)=g(0)+g'(v)u\ll\frac{1}{\theta^{2}}
\]
We apply this bound of $g(u)$ to the first integral of \eqref{eq:splitrange}
and obtain 
\[
\left|\int_{0}^{1/\sqrt{m}}e^{-(m+\alpha)u}u^{-\alpha}g(u)du\right|\ll\frac{1}{\theta^{2}}\int_{0}^{1/\sqrt{m}}e^{-(m+\alpha)u}u^{-\alpha}du.
\]
With substitution $(m+\alpha)u\rightarrow u$, the right side becomes
\begin{align*}
 & \frac{1}{\theta^{2}(m+\alpha)^{1-\alpha}}\int_{0}^{(m+\alpha)/\sqrt{m}}e^{-u}u^{-\alpha}du\\
\ll & \frac{\Gamma(1-\alpha)-\Gamma(1-\alpha,(m+\alpha)/\sqrt{m})}{\theta^{2}m^{1-\alpha}}
\end{align*}
where 
\[
\Gamma(s):=\int_{0}^{\infty}t^{s-1}e^{-t}dt\text{\ensuremath{\qquad\text{and}\qquad\Gamma(s,x):=\int_{x}^{\infty}t^{s-1}e^{-t}dt}}
\]
are the gamma and the upper incomplete gamma functions. With the asymptotic
formula for large $x\in\mathbb{R}^{+}$ and fixed $a\notin\mathbb{Z}$
\[
\Gamma(a,x)\sim x^{a-1}e^{-x}\sum_{n=0}^{\infty}\frac{\Gamma(a)}{\Gamma(a-n)}x^{-n},
\]
we obtain the following bound
\begin{equation}
\left|\int_{0}^{1/\sqrt{m}}e^{-(m+\alpha)u}u^{-\alpha}g(u)du\right|=\mathcal{O}\left(\frac{1}{\theta^{2}m^{1-\alpha}}\right).\label{eq:firstint}
\end{equation}
In summary, from \eqref{eq:splitrange}, \eqref{eq:secondint}, and
\eqref{eq:firstint}, we have an upper bound for the left side of
Lemma \ref{lem:keylemma} given by
\begin{equation}
\frac{1}{\pi}\int_{0}^{\infty}\frac{dx}{x^{\alpha}((1+x)^{2}-2(1+x)\cos\theta+1)(1+x)^{m+1}}\ll\frac{1}{\theta^{2}m^{1-\alpha}}.\label{eq:smallterm}
\end{equation}

To find a lower bound of the right side of this lemma, we note that
\[
\lim_{\theta\rightarrow0}\frac{\sin\theta(1-\cos\theta)^{\alpha}}{\sin^{\alpha}(\theta/2)}.\frac{1}{\theta^{1+\alpha}}=1
\]
and 
\[
\lim_{\theta\rightarrow\pi}\frac{\sin\theta(1-\cos\theta)^{\alpha}}{\sin^{\alpha}(\theta/2)}.\frac{1}{\theta^{1+\alpha}}=0.
\]
Thus the function 
\[
\frac{\sin\theta(1-\cos\theta)^{\alpha}}{\sin^{\alpha}(\theta/2)}.\frac{1}{\theta^{1+\alpha}}
\]
is bounded on $(0,\pi)$ and consequently 
\begin{equation}
\frac{\sin^{\alpha}(\theta/2)}{\sin\theta(1-\cos\theta)^{\alpha}}\gg\frac{1}{\theta^{1+\alpha}}.\label{eq:largeterm}
\end{equation}
By dividing each side of \eqref{eq:smallterm} by that of \eqref{eq:largeterm},
we conclude that 
\begin{equation}
\frac{\sin\theta(1-\cos\theta)^{\alpha}}{\sin^{\alpha}(\theta/2)}.\frac{1}{\pi}\int_{0}^{\infty}\frac{dx}{x^{\alpha}((1+x)^{2}-2(1+x)\cos\theta+1)(1+x)^{m+1}}\label{eq:fraction}
\end{equation}
is at most a constant multiple of 
\[
\frac{1}{(m\theta)^{1-\alpha}}.
\]
Thus there are constants $K$ and $M$ such that \eqref{eq:fraction}
is less than $1$ for all $\theta\in(K/m,\pi)$ and $m\ge M$ and
we complete the proof of Lemma \ref{lem:keylemma}.

\end{document}